\documentclass[11pt]{article}
\usepackage{latexsym,amsfonts,euscript}

\def\la{\langle}\def\ra{\rangle}

\def\pf{{\bf Proof.\quad }}
\def\pfend{\hfill{$\Box$}}

\title{A new graph related to conjugacy classes of finite groups
\thanks{Project supported by the NNSF of China (Grant No. 10571128),
the NSF of Hunan (Grant No. 04JJY4001) and the Scientific Research
Fund of Hunan Provincial Education Department (Grant No. 04C116).}}

\author{Xingzhong You \\
{\footnotesize\small School of Mathematics and Computing Science,
Changsha University }\\
{\footnotesize\small  of Science and Technology , Changsha, Hunan, 410077, P. R. China}\\
{\footnotesize\small  E-mail: xzyou2003@126.com}\\
\\
Guohua Qian\\
{\footnotesize\small Department of Mathematics,  Changshu
Institute of Technology,}\\
{\footnotesize\small  Changshu, Jiangsu, 215500, P.R.China}\\
{\footnotesize\small E-mail: ghqian2000@yahoo.com.cn} \\
\\
 Wujie Shi\\
 {\footnotesize\small School of Mathematics,  Suzhou
University,}\\
{\footnotesize\small Suzhou, Jiangsu, 215006, P.R. China}\\
{\footnotesize\small E-mail: wjshi@suda.edu.cn }}

\begin{document}
\date{}
\maketitle

\begin{abstract} In this paper we classify the finite groups satisfying the following property $P_4$:
their orders of representatives are set-wise relatively prime for
any 4 distinct non-central conjugacy classes.
\end{abstract}

\bigskip
{\bf 1. Introduction } Let $G$ be a finite group and $V$ the set of
the non-central conjugacy classes of $G$. From lengths of conjugacy
classes, the following class graph $\Gamma (G)'$ was introduced in
\cite{BHM}: its vertex set is the set $V$, and two distinct vertices
$x^G$ and $y^G$ are connected with an edge if $(|x^G|, |y^G|)>1$.
Similarly, in terms of orders of elements, we may attach a graph
$\Gamma (G)$ to $G$ as follows: its vertex set is also the set $V$,
and two distinct vertices $x^G$ and $y^G$ are connected with an edge
if $(o(x), o(y))>1$. Thus a new conjugacy class graph is defined.

The class graph $\Gamma (G)'$ has been studied in some detail: see
for example \cite{BHM}, \cite{CHM}, \cite{MP} and \cite{MQS}. In
\cite{MQS}, the authors have studied the structure of a finite group
$G$ with the following property: If for every prime integer $p$, $G$
has at most $n-1$ conjugacy classes whose size is a multiple of $p$.
In particular, they have classified the finite groups when $n=5$,
extending the result of Fang and Zhang \cite{MP}.

Inspired by \cite{MQS}, we study the structure of a finite group $G$
when its class graph $\Gamma(G)$ has no subgraph $K_n$, the complete
graph with $n$ vertices. First we observe from Lemma 1 that $\Gamma
(G)$ has no subgraph $K_n$ if and only if for any $n$ distinct
vertices $x_1^G, x_2^G, \dots, x_n^G$ in $\Gamma(G)$, the numbers
$o(x_1), o(x_2), \dots, o(x_n)$ are set-wise relatively prime, i.e.,
$(o(x_1), o(x_2), \dots, o(x_n))=1$. Consequently our problem can be
modified as follows. We say that a group $G$ satisfies property
$P_n$ if for every prime integer $p$, $G$ has at most $n-1$
non-central conjugacy classes whose order of representative is a
multiple of $p$. Thus $\Gamma (G)$ does not have a subgraph $K_n$ if
and only if $G$ satisfies property $P_n$.

The goal of this paper is to classify the finite groups that
satisfy property $P_4$.

\bigskip
{\bf Theorem A} \ \ \em Let $G$ be a finite group that satisfies
property $P_4$. Then $G$ is isomorphic to one of the following
groups:\em

\smallskip
(i)\ \ \em an abelian group;\em

(ii)\ \ \em The Frobenius group with complement of order 2 and
kernel $C_3$, $C_5$ or $C_7$;\em

(iii)\ \ \em The Frobenius group with complement of order 3 and
kernel $(C_2)^2$ or $C_7$;\em

(iv)\ \ \em $S_3\times C_2$ or $\la x, y|x^3=1, y^4=1,
xy=yx^{-1}\ra $;\em

(v)\ \ \em $D_8$ or $Q_8$;\em

(vi)\ \ \em The Frobenius group with complement of order 4 and
kernel $C_5$, $C_9$, $(C_3)^2$ or $C_{13}$; or\em

(vii)\ \ \em $L_2(5)$, $L_2(7)$, $L_2(9)$, $L_2(11)$, $L_2(13)$,
$A_7$, or $Sz(8)$.

\smallskip
Conversely, all these groups satisfy property $P_4$.\em

\bigskip
The proof of this theorem makes use of the classification theorem of
finite simple groups.

\bigskip

Throughout this paper, $\pi (n)$ denotes the set of all prime
divisors of a natural number $n$; $G$ denotes a finite group, set
$\pi (G)=\pi (|G|)$; for an element $x$ of $G$, $o(x)$ denotes the
order of $x$, and set $\pi (x)=\pi (o(x))$, $\pi_e (G)$ denotes the
set of all orders of elements in $G$; a class always means a
conjugacy class, and $x^G$ denotes the class containing $x$; if
$A\subseteq G$, let $k_G(A)$ be the number of classes of $G$ that
intersect $A$ nontrivially. All further unexplained notation is
standard.

\bigskip
{\bf 2. Preliminaries.} In this section, we will present some
preliminary results that we will need to prove Theorem A. We begin
with the following result.

\bigskip
{\bf Lemma 1.}  \em Let $G$ be a finite group. Then $G$ satisfies
property $P_n$ if and only if  $\Gamma(G)$ has no subgraph
$K_n$.\em

\bigskip
\pf  We only need to prove the necessity.

Suppose that $\Gamma (G)$ has a subgraph $K_n$ with $n$ vertices
$x_1^G, x_2^G, \dots, x_n^G$. If $Z(G)>1$, take $y(\not=1)\in
Z(G)$, then for $1\leq i\leq n$, we conclude that $x_iy\in G-Z(G)$
and $(x_iy)^G$ are $n$ distinct classes in $G$. But $o(y)|o(x_iy)$
for $1\leq i\leq n$, a contradiction. It follows that $Z(G)=1$.
Again, if $o(x_i)=p^k$ for some $i$ and prime $p$, then
$p|o(x_j)$, where $1\leq j\leq n$, also a contradiction.
Therefore, we have $|\pi(x_i)|\geq 2$ for every $i$. By induction
we can assume that $p|o(x_i)(1\leq i\leq n-1)$ for some prime $p$.
Let $z$ be the $p$-part of $x_1$. We conclude that $p|o(z)$ and
$z^G, x_1^G, x_2^G, \dots, x_{n-1}^G$ are $n$ distinct non-central
classes, also a contradiction.\pfend

\bigskip

{\bf Lemma 2.}  \em Let $G$ be a finite group that satisfies
property $P_n$. Then property $P_n$ is inherited by quotient
groups of $G$. \em

\bigskip
\pf Let $N$ be a normal subgroup of $G$ and $xN\in G/N$. Since
$o(xN)|o(x)$ and $(xN)^{G/N}$, when viewed as a subset of $G$, is
a union of some classes of $G$, the result follows.

\bigskip
We also use the following result that is taken from \cite{HuB}.

\bigskip
{\bf Lemma 3.} \em Suppose that $G$ is a non-abelian simple group
whose Sylow 2-subgroups are abelian. Then $G$ is isomorphic to one
of the following:

$L_2(q), q>3, q\equiv 3$, $5 (mod\, 8)$ or $q=2^f$, $J_1$, a Ree
group\em

\bigskip
From Lemma 3, we conclude that if $G$ is a non-abelian simple
group different from the ones in Lemma 3, then the Sylow
2-subgroups of $G$ are non-abelian and so $4\in \pi_e(G)$.

\bigskip
Recall the prime graph of a finite group $G$ is defined as
follows: its vertex set is $\pi (G)$, and two distinct vertices
$p$ and $q$ are connected with an edge if $pq\in \pi_e (G)$. Set
$\pi=\{p | 2$ is connected with $p$ and $p\in \pi(G) \}$.

\bigskip
In the following we determine the non-abelian simple groups that
satisfy property $P_4$.

\bigskip
{\bf Lemma 4.} \em Let $G$ be a non-abelian simple group that
satisfies property $P_4$. Then $G=$ $L_2(5)$, $L_2(7)$, $L_2(9)$,
$L_2(11)$, $L_2(13)$, $A_7$ or $Sz(8)$. \em

\bigskip
\pf  Since $G$ satisfies property $P_4$, we conclude from
\cite{Atlas} that $G$ is not of sporadic type. If $G$ is of
Alternating type, it follows from \cite{ZM} that $2, 4, 6, 8\in
\pi_e (G)$ provided that $n\geq 10$, a contradiction. Therefore,
$5\leq n\leq 9$ and hence $G=A_5, A_6$ or $A_7$ by \cite{Atlas}
again, where $A_5\cong L_2(5)$ and $A_6\cong L_2(9)$.

Now, assume that $G$ is of Lie type with characteristic $p$ and
$r$ is the number of classes of involutions in $G$. We split our
argument into two parts: $p$ is odd and $p=2$.

1.\ \ $p$ is odd.

(1)\ \ Let $G$ be of type $A_n(q), n\geq 1$.

We have $r\geq 2$ if $n\geq 3$. On the other hand, if $n\geq 2$,
then $4\in \pi_e (G)$. Also $|\pi|\geq 1$. It follows that $r=1$.
Therefore $n=1$ or $2$.

If $n=1$, then $G=A_1(q)\cong L_2(q)$, where $q\geq 5$. Let $x,
y\in G$ such that $o(x)=\frac{q-1}{2}, o(y)=\frac{q+1}{2}$, we
conclude that $|\pi(x)|\leq 2$ and $|\pi (y)|\leq 2$.

Suppose that $q\equiv 1(mod\, 4)$. Let $m=\frac{q-1}{4}$,
$s=\frac{q+1}{2}$. Since $(y^i)^G(1\leq i\leq m)$ are $m$ distinct
classes in $G$, we conclude that if $s$ is not prime, then $s\geq
9$ and hence $m\geq 4$, thus $G$ has at least four distinct
classes $y^G, (y^2)^G, (y^3)^G$ and $(y^4)^G$. Clearly
$o(y)=o(y^2)=o(y^4)=s$ and $1\neq o(y^3)|s$, a contradiction. This
implies that $s$ is prime and $s\leq 7$. It follows that $q=5, 9$
or $13$ and $G=L_2(5)$, $L_2(9)$ or $L_2(13)$.

Suppose that $q\equiv 3(mod\, 4)$. Let $m=\frac{q-3}{4}$,
$s=\frac{q-1}{2}$. Since $(x^i)^G(1\leq i\leq m)$ are $m$ distinct
classes in $G$, using a similar argument, we conclude that $q=7$
or $11$ and $G=L_2(7)$ or $L_2(11)$.

If $n=2$, then $G=A_2(q)$. It follows that there is a centralizer
of type $A_0(q)\times A_1(q)$. We have $|\pi|\leq 2$, a
contradiction.

(2)\ \ Let $G$ be one of the following types: $B_n(q)(n\geq 3)$,
$C_n(q)(n\geq 2)$, $D_n(q)(n\geq 4)$, $G_2(q)$, $F_4(q)$,
$E_6(q)$, $E_8(q)$, $^2A_n(q^2)(n\geq 2)$, $^2D_n(q^2)(n\geq 4)$
or $^2E_6(q^2)$.

Then $r\geq 2$ and $|\pi|\geq 1$. But $4\in \pi_e (G)$, a
contradiction.

(3)\ \ Let $G$ be of type $E_7(q)$.

Clearly, $4\in \pi_e(G)$. If $q=3$, then $r=3$, a contradiction. If
$q> 3$, then the prime graph of $G$ is connected. When $q\equiv 1
(mod\, 4)$, there is a maximal torus $T$ in $G$ such that
$|T|=\frac{1}{2}(q^3+1)(q^2+q+1)^2$. Clearly $\pi(|T|)\geq 3$. Since
$T$ is abelian, we can find an element $x\in T$ such that $|\pi
(x)|\geq 3$, a contradiction. When $q\equiv 3 (mod\, 4)$, there are
60 classes of maximal tori in $G$, all of which have even order,
then $|\pi|\geq 2$, also a contradiction.

(4)\ \ Let $G$ be of type $^3D_4(q^3)$.

Then $r=1$ and $|C_G(t)|= \frac{1}{2}q^4(q^2-1)(q^6-1)$ for an
involution $t$ in $G$. We have $|\pi|\geq 2$. Also, $4\in \pi_e
(G)$, a contradiction.

(5)\ \ Let $G$ be of type $^2G_2(q), q=3^{2m+1}, m\geq 1$.

By \cite[Chap XI Theorem 13.2]{HuB} , we have $r=1$ and
$C_G(t)=\la t \ra \times L_2(q)$ for an involution $t$ in $G$. It
follows that $|\pi|\geq 3$, a contradiction.

2.\ \ $p=2$.

(6)\ \ Let $G$ be of type $A_n(q), n\geq 1$.

Then $r=[\frac{n+1}{2}]$. Arguing as in (1), we have $r=1$, which
implies that $n=1$ or $2$.

If $n=1$, then $G=A_1(q)\cong L_2(q), q=2^f\geq 4$. Let $x, y\in
G$ such that $o(x)=q-1$ and $o(y)=q+1$. It follows that
$|\pi(x)|\leq 2$ and $|\pi(y)|\leq 2$. Since $(x^i)^G(1\leq i\leq
\frac{q-2}{2})$ are $\frac{q-2}{2}$ distinct classes in $G$,
arguing in a similar way in (1), we conclude that $q=4$ and
$G=L_2(4)$.

If $n=2$, take an involution $t\in G$, then $\pi(C_G(t))$
$=\pi(\frac{2(q-1)}{(3 ,\  q-1)})$. Note that there is a maximal
torus in $G$ of order $\frac{(q^2-1)}{(3 ,\  q-1)}$ and
$\pi(\frac{(q^2-1)}{(3 ,\  q-1)})\subseteq \pi_1$. If $q=2^f>4$,
then one of $q+1$ and $q-1$ is not prime. We conclude that there
is a prime divisor of $q-1$ such that it must be connected with
other three distinct primes in $\pi (G)$, a contradiction. It
follows that $q=2$ or $4$. By \cite{Atlas}, We have $G=L_3(2)\cong
L_2(7)$.

(7)\ \ Let $G$ be one of the following types: $C_n(q)(n\geq 2),
E_6(q)$ or $^2E_6(q^2)$.

Then $r=n+[\frac{n}{2}](\geq 3), 3$ or $3$, respectively. But
$4\in \pi_e (G)$, a contradiction.

(8)\ \ Let $G$ be one of the following types: $D_n(q)(n\geq 4)$,
$E_7(q)$, $E_8(q)$, $F_4(q)$ or $^2D_n(q^2)(n\geq 4)$.

Then $r=n+(-1)^n(\geq 4), 5, 4, 4$ or $2[\frac{n}{2}](\geq 4)$,
respectively, a contradiction.

(9)\ \ Let $G$ be one of the following types: $G_2(q), q> 2$,
$^3D_4(q^3)$ , $^2F_4(q), q>2$ or $^2F_4(2)'$.

Then $r=2$. But $4\in \pi_e (G)$ and $|\pi|\geq 1$, a
contradiction. \pfend

(10)\ \ Let $G$ be of type $^2A_n(q^2), n\geq 2$.

Then $r=[\frac{n+1}{2}]$. Since $4\in \pi_e (G)$ and $|\pi|\geq
1$, we conclude that $r=1$ and $n=2$. Let $t$ be an involution of
$G$. Then $\pi(C_G(t))=\pi(\frac{2(q^2+1)}{(3 ,\ q^2+1)})$. Note
that there is a maximal torus in $G$ of order $\frac{(q^4-1)}{(3
,\ q^2+1)}$ and $\pi(\frac{(q^4-1)}{(3 ,\ q^2+1)})\subseteq
\pi_1$. Arguing as in (6), we conclude that $q^2=4$ and $G=$
$^2A_2(4)$, a contradiction by \cite{Atlas}.

(11)\ \ Let $G$ be of type $^2B_2(q)(\cong Sz(q)), q=2^{2m+1},
m\geq 1$.

Applying \cite[Chap XI Lemma 11.6]{HuB} , we have that there are
$\frac{q}{2}-1$ classes of elements of order $q-1$ in $G$. If
$m\geq 3$, then $\frac{q}{2}-1\geq 63$, a contradiction. It
follows that $m=1$ or $2$ and hence $G=Sz(8)$ or $Sz(32)$. We have
$G=Sz(8)$ by \cite{Atlas}. \pfend

\bigskip
{\bf 3. Proof of Theorem A.} Now, we are ready to complete the
proof of Theorem A. It is easy to check that the groups listed in
Theorem A satisfy property $P_4$.

Now we split the classification into two parts. In Theorem 7 we
classify the groups that satisfy $P_4$ and $G'Z(G)<G$ and in
Theorem 8 we classify the ones that satisfy $P_4$ and $G'Z(G)=G$.

\bigskip
{\bf Theorem 5.} \em Let $G$ be a finite group with $G'Z(G)<G$. If
$G$ satisfies property $P_4$, then $G$ is one of the solvable
groups listed in \em(ii)-(vi) \em in Theorem A. \em

\bigskip
\pf Let $M=G'Z(G)$ and $p\in \pi(G/M)$. Take $xM\in G/M$ such that
$o(xM)=p$. Since $G/M$ is abelian, it follows that there are at
least $p-1$ classes of elements of order $p$ in $G/M$. Note that
$o(xM)|o(x)$ and $xM$, when viewed as a subset of $G$, is a union
of some classes of $G$, we conclude that $G$ has at least $p-1$
non-central classes whose order of representative is a multiple of
$p$. Therefore, $p-1\leq 3$, i.e., $p=2$ or $3$. Furthermore,
$|G/M|=2$, $3$ or $4$ and $k_G(G-M)\leq 3$.

1.\ \ Suppose that $k_G(G-M)=1$.

It follows from \cite[Proposition 2.1]{QSY} that $G$ is a
Frobenius group with kernel $M$ and $M$ is abelian of odd order
$\frac{|G|}{2}$. This implies that $Z(G)=1$ and $M=G'$. Since $G$
satisfies property $P_4$, we conclude that $M\in Syl_p(G)$ and
thus $k_G(M-\{1\})\leq 3$. It follows that $\frac{|M|-1}{2}\leq 3$
and hence $|M|\leq 7$. We deduce that $G$ is the group in (ii).

2.\ \ Suppose  that $k_G(G-M)=2$.

Applying \cite[Theorem 2.2]{QSY}, we get the following two cases.

(2.a)\ \ $|G/M|=3$ and $G$ is a Frobenius group with kernel $M$.

Similarly, we have $M\in Syl_p(G)$ and $k_G(M-\{1\})\leq 3$. If
$M$ is abelian, we conclude that $G$ is one of the groups in
(iii). If $M$ is non-abelian, then $k_G(Z(M)-\{1\})\leq 2$. Assume
first that $k_G(Z(M)-\{1\})=2$. We deduce that $|Z(M)|=7$ and $M$
is a $7$-group. Also we have that $M-Z(M)=x^G$ and $|x^G|=3\cdot
7^k$. Let $|M|=7^r$. Then $7^r=3\cdot 7^k + 7$, which has no
solution, a contradiction. Now, assume that $k_G(Z(M)-\{1\})=1$.
We have $|Z(M)|=4$ and $M$ is a $2$-group. Let $|M|=2^r$. If
$M-Z(M)=x^G$, then $|x^G|=3\cdot 2^k$ and hence $2^r=3\cdot 2^k +
4$, which implies $2^k=4$; if $M-Z(M)=x^G \cup y^G$, let
$|x^G|=3\cdot 2^k\leq 3\cdot 2^s=|y^G|$, then $2^r=3\cdot 2^s +
3\cdot 2^k+4$, which forces $(p^k, p^s, p^r)=(4,16, 64)$. We
conclude that there is always an element such that its centralizer
in $G$ is of order $4$. By \cite[Lemma 1.3]{QSY}, $M$ is the
dihedral, semi-dihedral or generalized quaternion group. This
forces $|Z(M)|=2$, also a contradiction.

(2.b)\ \ $|G/M|=2$ and $|C_G(x)|=4$ for any $x\in G-M$.

Applying \cite[Lemma 1.3 and Theorem 2.2(3)]{QSY}, we can see that
$Z(G)>1$. Since $|C_G(x)|=4$ for any $x\in G-M$, it follows that
$|Z(G)|=2$ . Take $x\in G-M$, we conclude that $o(xZ(G))=2$ and
$|C_{G/Z(G)}(xZ(G))|=2$ and thus $xZ(G)$ acts fixed point freely
on $M/Z(G)$, so $G/Z(G)$ is a Frobenius group with kernel
$M/Z(G)$. Clearly $M/Z(G)$ is a $p$-group, it follows that
$\frac{|M/Z(G)|-1}{2}\leq 3$ and hence $|M/Z(G)|=3,5$ or $7$, so
$|G|=12,20$ or $28$. We deduce that $G$ is the group in (iv)

3.\ \ Suppose that $k_G(G-M)=3$. Let $G-M=x^G\cup y^G\cup z^G$.

If $G$ is non-solvable, then \cite[Theorem 3.5]{QSY} will yield a
contradiction.

If $G$ is solvable, then by \cite[Theorem 3.6]{QSY} we get the
following three cases.

(3.a)\ \ $G\cong D_8$ or $Q_8$, thus $G$ is one of the groups in
(v).

(3.b)\ \ $G$ is a Froubenius group with kernel $M$ and cyclic
complement of order 4.

(3.c)\ \ $|G/M|=2$, $|C_G(x)|=|C_G(y)|=|C_G(z)|=6$, $o(x)=2$,
$o(y)=6$ and $z=y^{-1}$. And in this case, $M$ is of odd order and
$M$ has a normal and abelian 3-complement.

For case (3.b), arguing as in (1), we have $M\in Syl_p(G)$ and
$k_G(M-\{1\})\leq 3$. It follows that $\frac{|M|-1}{4}\leq 3$ and
hence $|M|\leq 13$. We conclude that $G$ is the group in (vi).

We claim that case (3.c) does not hold. If this is false, we can
see that $Z(G)=1$ and $M=G'$. Let $N$ be a normal and abelian
3-complement of $M$. Then $N$ is a normal and abelian
$\{2,3\}$-complement of $G$. If $N>1$, note that since there are 2
$G$-classes of elements of order 6 in $G-M$, then $M-N$ is exactly
one class of elements of order 3 in $G$. It follows from
\cite[Lemma 1.2]{QSY} that $M$ is a Frobenius group with kernel
$N$ and $|M/N|=3$. We conclude that $G/N\cong S_3$ and thus $G$ is
2-Froubenius. This forces $6\not \in \pi_e(G)$, a contradiction.
Hence $N=1$, thus $M$ is a 3-group and $k_G(M-\{1\})=1$.
Therefore, $M$ is abelian. This implies that for any $1\not=a\in
M$, we have $M\leq C_G(a)$. It follows that
$|M|-1=\frac{|G|}{|C_G(a)|}\leq |G/M|$, which forces that $|M|=3$
and so $|G|=6$, a contradiction. \pfend

\bigskip

{\bf Theorem 6} \em Let $G$ be a non-solvable finite group with
$G'Z(G)=G$. If $G$ satisfies property $P_4$, then $G$ is one of
the groups listed in \em (vii) \em in Theorem A. \em

\bigskip
\pf Since $G'Z(G)=G$,  it follows that
$(G/Z(G))'=G'Z(G)/Z(G)=G/Z(G)$. Thus there exists a normal
subgroup $N$ of $G$ with $N\geq Z(G)$ such that $G/N$ is a
non-abelian simple group. Since property $P_4$ is inherited by
quotient groups of $G$, we conclude from Lemma 4 that $G/N$ is
isomorphic to one of the following simple groups: $L_2(5)$,
$L_2(7)$, $L_2(9)$, $L_2(11)$, $L_2(13)$, $A_7$, $Sz(8)$.

Now it suffices to prove that $N$ is trivial. Suppose that this is
not true and let $G$ be a minimal counterexample. Then $N$ is a
minimal normal subgroup of $G$.

Suppose first that  $N$ is non-abelian. Then $N=N_1 \times \dots
\times N_s$ is a direct product of isomorphic simple groups $N_i$.
Observe that $Out(N)$ is solvable whenever $N$ is simple. We may
assume $s>1$. Take $x_1\in N_1$ and $y_1, y_2, y_3\in N_2$ such
that $o(x_1)=o(y_1)=2$, $o(y_2)=p$, $o(y_3)=q$, where $2,p$ and
$q$ are three distinct primes. We see that $x_1, x_1y_1, x_1y_2$
and $x_1y_3$ lie in distinct classes, a contradiction.

Suppose now that $N$ is a $p$-group and $C_G(N)>N$. Then $N=Z(G)$.
If $G'\cap N=1$, then $G=G'\times N$ since $G=G'Z(G)$, where $G'$
is the simple group listed above. By \cite{Atlas}, we can see
easily a contradiction. Therefore, $N=Z(G)\leq G'$ and so $G=G'$.
It follows that $N$ is the Schur multiplier of the simple group
$G/N$. Again by \cite{Atlas}, we can deduce a contradiction by
checking all possibilities of $G/N$. We conclude that $C_G(N)=N$.

Suppose then that $N$ is a $p$-group of order $p^r$ with $p\not\in
\pi(G/N)$ and $C_G(N)=N$. In this case, to find a contradiction we
take $G/N\cong L_2(5)$ as an example. Assume that $p \not\in
\pi(G/N)$. Then $G=NK$ by Schur-Zassenhaus' Theorem, where $K
\cong L_2(5)$. Let $Y<K$ be a Frobenius group order $12$, $x\in Y$
of order $3$, and $P$ a Sylow 2-subgroup of $Y$. Applying Theorem
15.16 of \cite{I},  we see that $|C_N(x)|=p^{r/3}$, and in
particular $3p\in \pi_e(G)$. Observe that $p, 2p\in \pi_e(G)$
since $P$ does not act fixed point freely on $N$. It follows that
$G$ has exactly one class of elements of order $3p$.  Clearly all
elements of order $3$ lie in one class. Thus $x^KN=x^G \cup y^G$,
where $y$ is of order $3p$. This implies that $$15p^r=
\frac{60p^r}{3p^{r/3}}+ \frac{60p^r}{|C_G(y)|},$$ \noindent and
thus $|C_G(y)|=12p^{r/3}/(3p^{r/3}-4)$, this is not an integer, a
contradiction.

Suppose finally  that $N$ is a $p$-group  with $p\in \pi(G/N)$ and
$C_G(N)=N$. If $G/N\cong L_2(5), L_2(7)$ or $L_2(9)$, then $|\pi
(G)|=3$ and so $k_G(G)\leq 10$; if $G/N\cong L_2(11), L_2(13)$ or
$A_7$, then $|\pi (G)|=4$, so $k_G(G)\leq 13$ and $k_G(G-N)\leq
10$. We conclude from \cite{VV} that $N=1$, a contradiction. Now
assume that  $G/N\cong Sz(8)$. Observe that $Sz(8)$ has  exactly
one class of involutions, two classes of elements of order $4$,
three classes of elements of order $7$, and three classes of
elements of order $13$, it follows that  $N$ is a $5$-group. Since
$G/N$ has Frobenius subgroups of order 14, we conclude that $10\in
\pi_e(G)$ by \cite[Lemma 1]{Maz}, and it follows that $G$ has at
least four distinct classes with representatives of even orders,
also a contradiction. \pfend

\bigskip

{\bf Remark.}\ \ In the proof of Theorem 6, we make use of the
classification of finite groups with few conjugacy classes given
in the papers in \cite{VV}. It would be possible to give an
independent proof, but this would make the proof longer.

\bigskip
\nopagebreak


\begin{thebibliography}{99}

\bibitem{BHM} Bertram, E. A., Herzog, M. and Mann, A., On a graph related to conjugacy classes of groups.
\em Bull. London Math. Soc. {\bf 22}, 569-575(1990).\em

\bibitem{CHM} Chillag, D., Herzog, M. and Mann, A., On the diameter of a graph related to conjugacy classes of groups.
\em Bull. London Math. Soc. {\bf 25}, 255-262(1993).\em

\bibitem{Atlas}   Conway, J. H., Curtis, R. I., Norton, S. P., Parker, R. A., Wilson, R. A.
``Atlas of finite groups". Oxford and New York: Oxford Univ.
 Press(Clarendon), 1985.

\bibitem{MP}  Fang, M., Zhang, P., Finite groups with graphs without triangles. \em J. Algebra {\bf 264}, 613-619(2003).\em

\bibitem{HuB} Huppert, B., Blackburn, N.,  ``Finite groups III". Berlin Heideberg New York: Springer-verlag, 1982.

\bibitem{I} Isaacs, I. M., Character theory of finite groups,  New
York: Academic Press, 1976.

\bibitem{VV} L\'{o}pez, A. V. and L\'{o}pez, J. V., Classification of finite groups according to the number of conjugacy
classes I, II. \em Israel. J. Math. {\bf 51}, 305-338(1985); {\bf
56}, 188-221(1986).\em

\bibitem{Maz} Mazurov, V. D., Characterizations of finite groups by sets of orders of their elements.
\em Algebra and Logik {\bf (1)36}, 23-32(1997).\em

\bibitem{MQS} Moret\'{o}, A., Qian, G. and Shi, W., Finite
groups whose conjugacy class graphs have few vertices. \em Arch.
math. {\bf 85}, 101-107(2005).\em

\bibitem{QSY} Qian, G.   Shi, W. and You, X., Conjugacy
classes outside a normal subgroup. \em Comm. in Algebra {\bf 32},
4809-4820(2004).\em

\bibitem{ZM} Zavarnitis, A. and Mazurov, V. D., Element orders in coverings of symmetric and alternating groups.
\em Algebra and Logik {\bf (3)38}, 159-170(1997).\em

\end{thebibliography}
\end{document}